\documentclass[12pt]{article}     \usepackage{amsmath}
\usepackage{amsthm}                \usepackage{amsopn}
\usepackage{graphics}
\usepackage{latexsym}              \usepackage{amsfonts}
\usepackage{amssymb}
\usepackage{amsmath}
\usepackage{enumerate}
\usepackage{array}
\usepackage{graphicx}
\usepackage{amsopn}

\swapnumbers
\theoremstyle{plain}
\newtheorem{theorem}{Theorem}[section]

\newtheorem{lemma}[theorem]{Lemma}
\newtheorem{corollary}[theorem]{Corollary}

\theoremstyle{definition}
\newtheorem{definition}[theorem]{Definition}
\newtheorem{remark}[theorem]{Remark}

\def\qed{\hfill\rule{1ex}{1ex}\\}
\newenvironment{pf}{\noindent {\bf Proof.}}{\qed}

\title{Reifenberg Flatness of Free\\
Boundaries in Obstacle Problems with\\
VMO Ingredients}
\author{Ivan Blank and Zheng Hao}

\newcommand{\nc}[2]{ \newcommand{#1}{#2} }

\nc{\avint}{ {- \hspace{-3.5mm} \int} }  
\nc{\xpr}{x^{\prime}}

\newcommand{\newavint}[1]{ \int_{#1} \mkern-37mu
                          \rule[.033 in]{.12 in}{.01 in} \ \ \  }
\newcommand{\myavinttwo}[1]{ \int_{#1} \mkern-35mu
                          \rule[.033 in]{.12 in}{.01 in} \ \ \  }
\nc{\R}{{\rm {I \! R}}}  
\nc{\N}{{\rm {I \! N}}}  
\newcommand{\closure}[1]{ \stackrel{\rule{.1 in}{.01 in}}{#1} }
\newcommand{\dclosure}[1]{ \stackrel{\rule{.2 in}{.01 in}}{#1} }

\newcommand{\chisub}[1]{ {\mathbf{\chi}}_{_{#1}} }

\newcommand{\newsec}[2]{ \section{#1} \label{sec-#2}  
                         \setcounter{equation}{0}     
                         \setcounter{theorem}{0} }    

\newcommand{\refeqn}[1]{ (\!\!~\ref{eq:#1}) } 
\newcommand{\refthm}[1]{ (\!\!~\ref{#1}) }    

\nc{\Holder}{H\"{o}lder\ }

\nc{\ith}{ \ensuremath{\text{i}^{\text{th}}} }
\nc{\jth}{ \ensuremath{\text{j}^{\text{th}}} }
\nc{\kth}{ \ensuremath{\text{k}^{\text{th}}} }
\nc{\curl}{ \nabla \times }
\nc{\Div}{ \nabla \cdot }

\nc{\Ppl}{ \mathcal{M}^{+} }  \nc{\Pmn}{ \mathcal{M}^{-} }

\nc{\smiley}{ $\stackrel{\because}{\smile} \;$ }

\newcommand{\BVP}[4]{  
  \begin{equation}
        \begin{array}{rl}
           #1 & \ \text{in}
               \ \ #4 \vspace{.05in} \\
           #2 & \ \text{on} \ \ \partial #4 \;.
        \end{array}
  \label{eq:#3}
  \end{equation}    }

\newcommand{\BVPb}[3]{   \BVP{#1}{#2}{#3}{ B_{1} }  }
\newcommand{\BVPhao}[4]{  
  \begin{equation}
        \left\{
        \begin{array}{ll}
           #1 & \ \text{in}
               \ \ #4 \vspace{.05in} \\
           \ \\
           #2 & \ \text{on} \ \ \partial #4 \;.
        \end{array} \right.
  \label{eq:#3}
  \end{equation}    }

\newcommand{\BVPc}[4]{  
  \begin{equation}
        \begin{array}{rl}
           #1 & \ \text{in}
               \ \ #4 \vspace{.05in} \\
           #2 & \ \text{on} \ \ \partial #4 \;,
        \end{array}
  \label{eq:#3}
  \end{equation}    }

\newcommand{\BVPcbsn}[3]{   \BVPc{#1}{#2}{#3}{ B_{s_0} }  }

\begin{document}
\numberwithin{equation}{section}
\maketitle

\begin{abstract} \noindent
We study the obstacle problem with an elliptic operator in divergence form.
We develop all of the basic theory of existence, uniqueness, optimal regularity, and nondegeneracy of the
solutions.  These results, in turn, allow us to begin the study of the regularity of the free boundary in the case
where the coefficients are in VMO.
\end{abstract}


\newsec{Introduction}{Intro}

\noindent
We study minimizers of
\begin{equation}
   \int_{B_1} a^{ij} D_i u D_j u
\end{equation}
among $u$ in the Hilbert space $W^{1,2}_0(B_1)$ which are constrained to lie above a fixed obstacle
$\varphi \in C^{0}(\closure{B_1}).$  (We use Einstein summation notation throughout the paper.)
We assume that our obstacle $\varphi < 0$ on $\partial B_1,$ and to avoid triviality we will assume
that $\max \varphi > 0.$
We assume that at each $x \in B_1,$ the matrix $\mathcal{A} = (a^{ij})$ is symmetric and strictly and uniformly elliptic, i.e.
\begin{equation}
  \mathcal{A} \equiv \mathcal{A}^{T} \ \ \text{and} \ \ 0 < \lambda I \leq \mathcal{A} \leq \Lambda I \;,
\label{eq:UniformEllip}
\end{equation}
or, in coordinates:
$$a^{ij} \equiv a^{ji} \ \ \text{and} \ \
    0 < \lambda |\xi|^2 \leq a^{ij} \xi_i \xi_j \leq \Lambda |\xi|^2 \ \ \text{for all} \ \xi \in \R^n, \ \xi \ne 0 \;.$$

If we let $Lv := D_{i}a^{ij}D_{j} v$ in the usual weak sense for a divergence form operator and we
consider the case where $L\varphi \in L^{\infty}(B_1),$ then by letting $w:= u - \varphi$ and by letting
$f := -L\varphi,$ the study of the minimizers above leads us to look at
weak solutions of the obstacle-type problem:
\begin{equation}
           Lw := D_{i}a^{ij}D_{j}w = \chisub{ \{ w > 0 \} }f \ \ \text{in}
               \ \ B_1 \;,
\label{eq:BasicProb}
\end{equation}
where $\chisub{S}$ denotes the characteristic function of the set $S,$ and where we look for
$w \geq 0.$ A weak solution to a
second order partial differential equation is a weakly differentiable function which satisfies
an appropriate equality when integrated against test functions.  (See chapter 8 of \cite{GT}.)
As an example, we will say that $w \in W^{1,2}(B_1)$ satisfies Equation\refeqn{BasicProb}if
for any $\phi \in W^{1,2}_{0}(B_1)$ we have:
\begin{equation}
           - \int_{B_1} a^{ij}D_{j}w D_{i}\phi = \int_{B_1} \phi \chisub{ \{ w > 0 \} }f  \;.
\label{eq:BasicProbWeakForm}
\end{equation}

Our motivations for studying this type of problem are primarily theoretical.  Indeed,
the obstacle problem is possibly the most fundamental and important free boundary problem,
and it originally motivated the study of variational inequalities.  On the other hand, the obstacle problem
has well-established connections to the Stefan problem and the Hele-Shaw problem.  (See \cite{C1} and
\cite{BKM} for example.)  Furthermore, as observed in \cite{MPS} the mathematical modeling of numerous
physical and engineering phenomena can lead to elliptic problems with discontinuous coefficients, and so
the current case seems to allow some of the weakest possible solutions.

Our main result is the following:
\begin{theorem}[Free Boundary Regularity]  \label{FBRi}
We assume
\begin{enumerate}
    \item $w \geq 0$ satisfies Equation\refeqn{BasicProb}\!\!,
    \item $a^{ij}$ satisfies Equation\refeqn{UniformEllip}\!\!,
    \item $0 < \lambda^{\ast} \leq f \leq \Lambda^{\ast},$ and
    \item $a^{ij}$ and $f$ belong to the space of vanishing mean oscillation (VMO).
\end{enumerate}
We let $S_r$ denote the set of regular points of the free boundary within $B_r,$ and
assume $K \subset \subset S_{1/2}.$  Then $K$ is a Reifenberg vanishing set.
\end{theorem}
\noindent
The definition of Reifenberg vanishing is found at the beginning of the fifth section.

As a corollary of this result we will conclude that blowup limits at regular points will be rotations and
scalings of the function $(x_n^{+})^2.$  In terms of the fact that this function is homogeneous of degree
2, it is quite usual to use Weiss's celebrated monotonicity formula to prove this type of result. (See \cite{W}.)
On the other hand, the weak nature of our equation, together with the weak $W^{1,2}$ convergence
to blowup solutions make it difficult to estimate differences of the values of the Dirichlet integrals which
appear in Weiss's formula.  So, instead of using homogeneity to prove Reifenberg flatness, our paper goes
in the opposite direction.


\newsec{Preliminaries and Basic Results}{PrelBasRes}

\noindent
We will use the following basic notation throughout the paper:
$$
\begin{array}{lll}
\chisub{D} & \ & \text{the characteristic function of the set} \ D \\
\closure{D} & \ & \text{the closure of the set} \ D \\
\partial D & \ & \text{the boundary of the set} \ D \\
x   & \ & (x_1, x_2, \ldots, x_n) \\
x^{\prime} & \ &(x_1, x_2, \ldots, x_{n-1}, 0) \\
B_{r}(x) & \ & \text{the open ball with radius} \ r \ \text{centered at the point} \ x \\
B_{r}  & \ & B_{r}(0) \\
\Omega(w) & \ & \{ w > 0 \} \\
\Lambda(w) & \ & \{ w = 0 \} \\
FB(w) & \ & \partial \Omega(w) \cap \partial \Lambda(w) \\
\end{array}
$$
Throughout the entire paper, $n, \lambda,$ and $\Lambda$ will remain fixed, and so
we will omit all dependence on these constants in the statements of our theorems.
We will typically work in the Sobolev spaces and the \Holder spaces, and we will follow all of the definitions and
conventions found in the book by Gilbarg and Trudinger.  (See \cite{GT}.)  To simplify
exposition slightly, for $u,v \in W^{1,2}(D)$ we will say that $u = v$ on $\partial D$ if
$u - v \in W^{1,2}_{0}(D).$

We define the divergence form elliptic operator
\begin{equation}
     L := D_j \; a^{ij}(x) D_i \;,
\label{eq:Ldef}
\end{equation}
or, in other words, for a function $u \in W^{1,2}(\Omega)$ and $f \in L^2(\Omega)$ we say ``$Lu = f$ in $\Omega$'' if
for any $\phi \in W_{0}^{1,2}(\Omega)$ we have:
\begin{equation}
    - \int_{\Omega} a^{ij}(x) D_{i} u D_{j} \phi = \int_{\Omega} g \phi \;.
\label{eq:Ldef2}
\end{equation}
(Notice that with our sign conventions we can have $L = \Delta$ but not $L = -\Delta.$)
Next, we fix a function $\psi \in W_{loc}^{1,2}(\R^n)$ with $\psi \geq 0$ which we will use as boundary data,
and we fix a function $\varphi \in C^{0}(\dclosure{B_1})$ which we will use as an obstacle.

%

Define the functionals:
$$D(u, \Omega) := \int_{\Omega} (a^{ij}D_i u D_j u) \;, \ \ \text{and}$$
$$J(w, \Omega) := \int_{\Omega} (a^{ij}D_i wD_j w + 2w) \;.$$
For any bounded set $\Omega \subset \R^n$ we will minimize these functionals in the following sets, respectively:
$$S_{\Omega, \varphi} := \{ u \in W^{1,2}_{0}(\Omega) \; : u \geq \varphi \; \} \;, $$
$$H_{\Omega,\psi} := \{ w \in W^{1,2}(\Omega) \; : \; w - \psi \in W_{0}^{1,2}(\Omega) \; \} \;, \ \ \text{and}$$
$$K_{\Omega, \psi} := \{ \ w \in H_{\Omega, \psi} \; : \; w(x) \geq 0 \ \text{for all} \ x \in \Omega \; \}.$$
When it is clear on which set we are working, we will simply write ``$D(u)$'' in place of ``$D(u, \Omega)$'' and
``$S_{\varphi}$'' in place of ``$S_{\Omega, \varphi}$'' and so on. 

Probably the most classic version of the obstacle problem involves minimizing $D(u, B_1)$ within $S_{B_1,\varphi}$ in
the case where $a^{ij} = \delta^{ij}.$  (Here we use $\delta^{ij}$ to denote the usual Kronecker delta function so that
$D(u)$ simplifies to the usual Dirichlet integral.  See \cite{C1}, \cite{C2}, \cite{C4}, and \cite{C5} for an analysis of this
problem.)  Indeed, following the same arguments given at the beginning of \cite{C5}, but for the more general $a^{ij}$
considered here, we can establish the following theorem:
\begin{theorem}[Basic Results]  \label{BasicRes}  Given an obstacle $\varphi \in W^{1,2}(B_1)$ which has a trace on
$\partial B_1$ which is negative almost everywhere,
there is a unique $u \in S_{B_1,\varphi}$ which minimizes $D(u, B_1).$  Furthermore, $u$ is a bounded supersolution
to the problem $L(u) = 0.$ Finally, if $\varphi$ is continuous, then $u$ is almost everywhere equal to a function which is
continuous on all of $\closure{B_1}.$
\end{theorem}

\begin{pf}  For the proof, just follow the beginning of \cite{C5}.  (Note that the details of the proof of the mean value
formula that Caffarelli uses can be found within \cite{BH}.)
\end{pf}


Turning to the regularity questions, we find it convenient to work with the height function $w$ which is the minimizer
of $J$ within $K_{B_1,\psi}.$  On the other hand, one can ask if this is really the same problem as before.  In the original
problem with the Laplacian (in other words, with $a^{ij} = \delta^{ij}$), if the obstacle is twice differentiable, then it
makes sense to take its Laplacian.  In the current situation, it is not as simple to characterize the functions $\varphi,$
where $L\varphi$ makes sense.  The obvious route, however, is to simply assume that $L \varphi = -f$ for a function
$f$ with specified properties.  If we assume that $L\varphi = -f,$ and that $f \in L^{\infty}(B_1),$ then the two problems
are completely equivalent.

We are most interested in the obstacle problem where we minimize $J$ within $K_{B_1,\psi}.$  Besides requiring existence
and regularity, we need to know that the minimizer, $w,$ satisfies $w \geq 0$ and
\BVPb{L(w) = \chisub{ \{ w > 0 \} }f}{w = \psi}{OPPDE}
The proof of this fact and many of the related facts follows \cite{BH} very closely, and so we will only mention that
the proof is carried out with a penalization argument.  The details can be found with only very minor adjustments in
\cite{BH}.  To summarize the relevant facts we can state the following result:

\begin{theorem}[Problem Equivalencies] \label{ProbEquiv}
Let $\varphi$ be an obstacle which satisfies the following:
\begin{enumerate}
    \item $\psi := -\varphi > 0$ on all of $\partial B_1.$
    \item $f := -L\varphi \in L^{\infty}(B_1).$
\end{enumerate}
Finally assume that $w = u - \varphi.$  Then the following are equivalent:
\begin{enumerate}
    \item $w$ satisfies Equation\refeqn{OPPDE}\!\!.
    \item $w$ minimizes $J$ in $K_{B_1,\psi}.$
    \item $u \in W^{1,2}_{0}(B_1)$ satisfies $Lu = -\chisub{ \{ u = \varphi \} }f.$
    \item $u$ minimizes $D$ in $S_{B_1,\psi}.$
\end{enumerate}
\end{theorem}


\noindent
Now in order to get to the regularity of the free boundary we need two more basic facts which can also be found
within \cite{BH}.  At this point, having proven our theorem about the equivalencies between the problems, it is worth
gathering a collection of assumptions that we will have for the rest of this paper.  We will always assume:
\begin{equation}
   \begin{array}{l}
       \displaystyle{L(w) = \chisub{ \{ w > 0 \} }f \ \ \text{in} \ B_1 \;,} \\
\ \\
       \displaystyle{a^{ij}(x) \equiv a^{ji}(x) \;,} \\
\ \\
       \displaystyle{
    0 < \lambda |\xi|^2 \leq a^{ij} \xi_i \xi_j \leq \Lambda |\xi|^2 \ \ \text{for all} \ \xi \ne 0 \;, } \\
\ \\
       \displaystyle{0 < \lambda^{\ast} \leq f \leq \Lambda^{\ast} \;, \ \ \text{and} } \\
\ \\
       \displaystyle{w \geq 0}
   \end{array}
\label{eq:Always}
\end{equation}
and we will frequently assume
\begin{equation}
   0 \in \partial \{w > 0\}.
\label{eq:Freq}
\end{equation}

For the next two theorems we assume both Equation\refeqn{Always}and Equation\refeqn{Freq}\!\!.
We start with a regularity statement which gives us compactness of quadratic rescalings.


\begin{theorem}[Optimal Regularity]  \label{opreg}
For any $x \in B_{1/2}$ we have
\begin{equation}
   w(x) \leq \tilde{C}|x|^2
\label{eq:OpReg}
\end{equation}
for a constant $\tilde{C} = \tilde{C}(n, \lambda, \Lambda, \lambda^{\ast}, \Lambda^{\ast} ).$
\end{theorem}

On the other hand, there is a nondegeneracy statement which prevents quadratic rescalings from vanishing
in the blow up limit.  Namely, we have:


\begin{theorem}[Nondegeneracy] \label{NonDeg}
With $C = C(n,\lambda, \Lambda, \lambda^{\ast}, \Lambda^{\ast}) > 0,$ and for any $r \leq 1$ we have
\begin{equation}
    \sup_{x \in B_r} w(x) \geq Cr^2 \;.
\label{eq:NonDeg}
\end{equation}
\end{theorem}

Although the optimal regularity statement can be proven by a straightforward adjustment of the proof for the case
when $a^{ij} \equiv \delta^{ij},$ the proof of nondegeneracy is much easier in the case with the Laplacian because
of the usefulness of the function $|x|^2.$  In the present case, in order to prove nondegeneracy one seems to
need a polygonal curve argument and this can be found in \cite{BH}.


\newsec{Measure Stability}{meastab}

Now we begin a measure theoretic study of regularity which will culminate in a measure theoretic version
of the theorem proven by Caffarelli in 1977.  (See \cite{C1}.)

\begin{lemma}[Compactness I] \label{Lpoint}
Let $\{a^{ij}_k\},$ $\{f_k\},$ and $\{w_k\}$ satisfy
\begin{enumerate}
  \item $0 < \lambda I \leq a^{ij}_k \leq \Lambda I,$
  \item $0 < \lambda^{\ast} \leq f_k \leq \Lambda^{\ast},$
  \item $w_k \geq 0$, $D_ia^{ij}_kD_jw_k=\chisub{ \{ w_k>0 \} }f_k$ in $B_2,$ and $0\in\partial\{w_k>0\}$.
  \item $||w_k||_{W^{1,2}(B_2)} \leq \gamma < \infty$ and 
  \item there exists an $f$ (with $0 < \lambda^{\ast} \leq f \leq \Lambda^{\ast}$),
           such that $f_k$ converges to $f$ strongly in $L^1.$
\end{enumerate}
then there exists a $w \in W^{1,2}(B_1)$ and an $f \in L^{\infty}(B_1)$ and a
subsequence of $\{w_k\}$ such that along this subsequence (which we still label with ``k''), we have
\begin{itemize}
   \item[A.] uniform convergence of $w_k$ to $w,$ and weak convergence in $W^{1,2},$
   \item[B.] for any $\phi\in W^{1,2}_0(B_1)$
\begin{equation}
\label{eq:goodchiconv}
   \int_{B_1}\chisub{ \{ w_k>0 \} } f_k \phi \rightarrow \int_{B_1}\chisub{ \{ w>0 \} } f \phi.
\end{equation}
\end{itemize}
\end{lemma}

\begin{pf}  Item A follows by using standard functional analysis combined with De Giorgi-Nash-Moser theory.
Since we can take a subsequence, we can assume without loss of generality that $f_k$ converges to $f$
pointwise almost everywhere.
In the interior of both $\{ w > 0 \}$ and $\{ w = 0 \}$ it is not hard to show that $\chisub{ \{ w_k>0 \} } f_k$
converges pointwise almost everywhere to $\chisub{ \{ w>0 \} } f$
(for the interior of $\{ w = 0 \}$ one needs to use the nondegeneracy statement),
so by Lebesgue's dominated convergence theorem it suffices to prove that $\partial \{ w = 0 \}$ has no Lebesgue points.
The proof of this fact is very similar to the proof of Lemma 5.1 of \cite{BT}, but we include it here for the convenience of
the reader.

Let $x_0\in\partial\{w=0\} \cap B_1,$ and choose $r > 0$ such that
$$B_r(x_0)\subset B_1.$$
Define $W(x):= r^{-2}w(x_0+rx)$ and $W_k(x):= r^{-2}w_k(x_0+rx)$. After this change of
coordinates, we have $0\in\partial\{W=0\},$ and so there exists
$\{x_k\}\rightarrow 0$ such that $$W(x_k)>0, \ \text{for all}\ k.$$
Now fix $k$ so $x_k\in B_{1/8}$, take $J$ large enough such that $i,j\geq J$ implies
\begin{equation}
\label{eq:WjvsW}
  ||W_j-W||_{L^{\infty}(B_1)}\leq \frac{W(x_k)}{2},
\end{equation}
and
\begin{equation}
\label{eq:WivsWj}
  ||W_i-W_j||_{L^{\infty}(B_1)}\leq \frac{\tilde{C}}{10}
\end{equation}
where $\tilde{C}= \frac{C}{10}$ which is the constant from the nondegeneracy statement.

Since $W_j\rightarrow W$ in $C^{\alpha}$, $W_J(x_k)>0$ and nondegeneracy imply the existence
of $\tilde{x}\in B_{1/2}$ such that
\begin{equation}
\label{eq:towrdcontra}
   W_J(\tilde{x}) \geq C \left(\frac{1}{2} - \frac{1}{8} \right)^2 = \frac{9}{64} C > \tilde{C}.
\end{equation}
Now $i\geq J$ implies $W_i(\tilde{x})\geq \frac{9\tilde{C}}{10}$.
Since $W_i$ satisfies a uniform $C^{\alpha}$ estimate, there exists an $\tilde{r}>0$ such that
$W_i(y)\geq \frac{\tilde{C}}{2}$ for all $y\in B_{\tilde{r}}(\tilde{x})$ once $i\geq J$.
From this we can conclude $B_{\tilde{r}}(\tilde{x})\subset \{W_{\infty}>0 \}$.

Scaling back to the original functions, we conclude $x_0$ is not Lebesgue point.
Since $x_0$ was an arbitrary point of the free boundary there are no Lebesgue points
in $\partial\{w>0\}$.
\end{pf}


\begin{lemma}[Compactness II] \label{CompII}
If we assume everything we did in the previous lemma, and we assume in addition that
$A = (A^{ij})$ is a symmetric, constant matrix with
$$0 < \lambda I \leq A \leq \Lambda I,$$
and such that
$$||a^{ij}_k - A ^{ij}||_{L^1(B_1)} \rightarrow 0,$$
then the limiting functions $w$ and $f$ given in the last lemma satisfy:
\begin{equation}
     D_i A^{ij} D_j w = \chisub{ \{ w > 0 \} }f
\label{eq:FinalEqn}
\end{equation}
in $B_1.$  Furthermore, $0 \in \partial \{ w > 0 \}.$
\end{lemma}

\begin{pf}
Since $a^{ij}_k\rightarrow A^{ij}$, and there is a uniform $L^{\infty}$ bound on all of $a^{ij}_k$ and $A^{ij}$, we have
\begin{equation}
\label{eq:strgLq}
a^{ij}_k\rightarrow A^{ij} \ \ \text{in}\ L^{q}(B_1)
\end{equation}
for any $q<\infty$, in particular $a^{ij}_k\rightarrow A^{ij}$ in $L^2$. We have for any $\phi\in W^{1,2}_{0}(B_1)$,
\begin{align*}
\int_{B_1}a^{ij}_kD_i w_k D_j\phi & =\int_{B_1}(a^{ij}_k-A^{ij})(D_i w_k-D_i w)D_j\phi\\
                                  & +\int_{B_1}a^{ij}_k D_i w D_j\phi + \int_{B_1} A^{ij}(D_iw_k-D_iw)D_j\phi
\end{align*}
Since $a^{ij}_k\rightarrow A^{ij}$ in $L^2$ and $D_iw_k \rightharpoonup D_iw$, we have
\begin{equation} \int_{B_1}(a^{ij}_k-A^{ij})(D_i w_k-D_i w)D_j\phi \rightarrow 0,
\end{equation}
\begin{equation} \int_{B_1}a^{ij}_k D_i w D_j\phi \rightarrow \int_{B_1}A^{ij}D_iwD_j\phi
\end{equation}
and
\begin{equation}
\int_{B_1} A^{ij}(D_iw_k-D_iw)D_j\phi \rightarrow 0.
\end{equation}
Therefore,
\begin{equation}
\int_{B_1}a^{ij}_kD_i w_k D_j\phi \rightarrow \int_{B_1}A^{ij}D_i w D_j\phi .
\end{equation}
Together with Equation\refeqn{goodchiconv}\!\!, we proved 
$$D_i A^{ij} D_j w = \chisub{ \{ w > 0 \} }f.$$

Now in order to show that $0 \in \partial \{w > 0\}$ we observe first that
$0\in\partial\{w_k > 0\}$ implies $$0\in\{w=0\}.$$
Next we suppose there exists $r_0$, such that $B_{2r_0}\in\{w=0\}$. For any $k$, we have 
\begin{equation}
    \sup_{x \in B_{r_0}} w_k(x) \geq C(r_0)^2 \;.
\label{eq:NonDeg}
\end{equation}
By picking a convergent subsequence we get a contradiction to
$w=0$ in $B_{2r_0}.$ Therefore, we have $0 \in \partial \{ w > 0 \}.$
\end{pf}


\begin{theorem}[Measure Stability] \label{MeaStab}
Fix positive constants $\gamma, \lambda, \Lambda, \lambda^{\ast},$ and $\Lambda^{\ast},$
and suppose $w$ satisfies Equation\refeqn{Always}\!\!,
and for some constant $\mu \in [\lambda^{\ast}, \Lambda^{\ast}],$ assume that $u$ satisfies
\begin{equation} \Delta u =\chisub{ \{ u>0 \} } \mu \ \ \text{in} \ B_1 \;
\end{equation}
with $$w=u, \ \ \text{on} \ \partial B_1.$$
where we assume in addition that $w$ satisfies $$||w||_{W^{1,2}(B_1)}\leq \gamma, \ \ \text{and} \ \
||w||_{C^{\alpha}(\closure{B_1})} \leq \gamma.$$

Then there exists a modulus of continuity $\sigma(\epsilon)$, such that if
\begin{equation}
     ||a^{ij}-\delta^{ij}||_{L^2(B_1)} < \sigma(\epsilon), \ \ \ \text{and} \ \ \
     ||f - \mu||_{L^1(B_1)} < \sigma(\epsilon)
\end{equation}
 then
\begin{equation}
   |\{w=0\}\Delta\{u=0\}| < \epsilon.
\end{equation}
(We are abusing notation slightly by using $\mu$ to denote the function which is everywhere equal to $\mu$
in $B_1.$)
\end{theorem}

\begin{pf} The proof of Theorem 5.4 of \cite{BT} can be adapted to the current setting without too much
difficulty, but we include it for the convenience of the reader.
Suppose not. Then there exists $a^{ij}_k, \; w_k, \; f_k$ and $u_k$ such that,
 \begin{enumerate}
  \item $D_i a^{ij}_k D_j w_k=\chisub{ \{ w_k>0 \} } f_k$ in $B_1,$
  \item $a^{ij}_k\rightarrow \delta^{ij}$ in $L^2(B_1),$
  \item $f_k \rightarrow \mu$ in $L^1(B_1),$
  \item $\begin{cases}
                  \Delta u_k=\chisub{ \{ u_k>0 \} } \mu \ & \text{in}\  B_1 \\
                  u_k=w_k \ & \text{on} \ \partial B_1, \ \ \ \text{and}
             \end{cases}$
  \item $||w_k||_{W^{1,2}(B_1)}\leq \gamma,$ and $||w_k||_{C^{\alpha}(\closure{B_1})} \leq \gamma.$
 \end{enumerate}
but $ |\{w_k=0\}\Delta\{u_k=0\}| \geq \epsilon_0 $ for some $\epsilon_0$ fixed.

By applying the previous compactness lemmas to an arbitrary subsequence, there exists a $w_{\infty}$
and a sub-subsequence such
that $$w_k \rightharpoonup w_{\infty}, \ \ \text{in} \ \ W^{1,2}(B_1)$$ and
$$w_k\rightarrow w_{\infty} \ \ \text{in} \ \ C^{0}(\closure{B_1})$$
which implies $w_k \rightarrow w_{\infty}$ in $L^2(B_1).$  (We will still use ``$w_k$'' for the sub-subsequence.)
Equation\refeqn{goodchiconv}is also satisfied with the constant function $\mu$ in place of $f.$

By standard comparison results for the obstacle problem (see for example Theorem 2.7a of \cite{B}),
there exists $u$ such that
\begin{equation}
\label{eq:ukconvtou}
    u_k \rightarrow u \ \text{in} \ L^{\infty}(B_1)
\end{equation}
We have for any $\phi\in W^{1,2}_{0}(B_1)$,
\begin{align*}
\int_{B_1} a^{ij}_k D_i w_k D_j\phi & =\int_{B_1}(a^{ij}_k - \delta^{ij})( D_i w_k - D_i w_{\infty} )D_j \phi\\
        & + \int_{B_1} a^{ij}_k D_i w_{\infty} D_j \phi + \int_{B_1} \delta^{ij} (D_i w_k - D_i w_{\infty}) D_j \phi
\end{align*}
Since $a^{ij}_k \rightarrow \delta^{ij}$ in $L^2$ and $D_i w_k \rightharpoonup D_i w_{\infty}$, we have
\begin{equation} \int_{B_1} (a^{ij}_k - \delta^{ij})(D_i w_k - D_i w_{\infty}) D_j\phi \rightarrow 0,
\end{equation}
\begin{equation} \int_{B_1} a^{ij}_k D_i w_{\infty} D_j \phi \rightarrow \int_{B_1} \delta^{ij} D_i w_{\infty} D_j \phi
\end{equation}
and
\begin{equation}
\int_{B_1} \delta^{ij}( D_i w_k - D_i w_{\infty} )D_j\phi \rightarrow 0.
\end{equation}
Therefore,
\begin{equation}
\int_{B_1}a^{ij}_kD_i w_k D_j\phi \rightarrow \int_{B_1}\delta^{ij} D_i w_{\infty} D_j\phi .
\end{equation}
By Equation\refeqn{goodchiconv}with $\mu$ in place of $f,$ we have
\begin{equation}
   \int_{B_1} \chisub{ \{ w_k>0 \} } f_k \phi  \rightarrow \int_{B_1} \chisub{ \{ w_{\infty} > 0 \} } \mu \phi,
\end{equation}
so $w_{\infty}$ satisfies
\begin{equation}
\Delta w_{\infty} =\chisub{ \{ w_{\infty} > 0 \} } \mu \ \ \text{in} \ B_1.
\end{equation}
We notice that by assumption,
\begin{align*} 0<\epsilon_0 & \leq |\{w_k=0\}\Delta\{u_k=0\}| \\
                            & = ||\chisub{ \{ u_k>0 \} }-\chisub{ \{ w_k>0 \} }||_{L^1(B_1)} \\
                            & \leq ||\chisub{ \{ u_k>0 \} }-\chisub{ \{ w_{\infty}>0 \} }||_{L^1(B_1)} +
                                      ||\chisub{ \{ w_{\infty}>0 \} }-\chisub{ \{ w_k>0 \} }||_{L^1(B_1)} \\
                            &= I+I\!I \;.
\end{align*}
For $I$, since \BVPhao{\Delta u_k = \chisub{ \{u_k > 0\} } \mu}{u_k = w_k}{PDEuk}{B_1}
and \BVPhao{\Delta w_{\infty} = \chisub{ \{ w_{\infty} > 0\} }\mu}{w_{\infty} = u}{PDEw}{B_1}
By Theorem 2.7a of \cite{B},
we have
\begin{equation}
    ||u_k - w_{\infty} ||_{L^\infty(B_1)}\leq ||u_k-u||_{L^\infty(\partial B_1)} \;,
\label{eq:512analog}
\end{equation}
and since $u_k\rightarrow u$ in $L^\infty,$
we have 
\begin{equation}
   ||\chisub{ \{ u_k>0 \} } - \chisub{ \{ w_{\infty}>0 \} }||_{L^1(B_1)} \rightarrow 0,
\label{eq:CHIniceconv}
\end{equation}
by Corollary 4 of \cite{C2}.

For $I\!I,$ we know that inside $\{ w_{\infty}>0 \},$  $w_k$ will eventually be positive by the uniform convergence,
so $\chisub{ \{ u_k>0 \} }= \chisub{ \{ w_{\infty} > 0 \} }$ there.
In the interior of $\{ w_{\infty} = 0 \},$  $w_k$ will eventually be
$0,$ since otherwise we will violate the nondegeneracy property, and so
$\chisub{ \{ u_k > 0 \} }= \chisub{ \{ w_{\infty} > 0 \} }$ there.
Finally, since $\partial \{w_{\infty}=0\}$ has finite $(n-1)$-dimensional Hausdorff measure (see \cite{C2}, \cite{C3},
or \cite{C5}), we must have $|\partial \{w_{\infty}=0\}| = 0,$ and therefore $I\!I \rightarrow 0.$
This convergence to $0$ gives us a contradiction, since $0 < \epsilon_0 \leq I + I\!I.$
\end{pf}




\newsec{Weak Regularity of the Free Boundary}{WRFB}
In this section we establish the existence of blow up limits, and use this result to show a measure-theoretic version
of Caffarelli's free boundary regularity theorem. We will show the existence of blowup limits in the case where the
$a^{ij}$ and the $f$ belong to VMO.  We define VMO to be the subspace of BMO such that if $g \in BMO$ and
\begin{equation}
    \eta_{g}(r) : = \sup_{\rho \leq r, \; y \in \R^n} \ \frac{1}{|B_{\rho}|} \int_{B_{\rho}(y)} |g(x) - g_{_{B_{\rho}(y)}}| \; dx \;,
\label{eq:etadef}
\end{equation}
then $\eta_{g}(r) \rightarrow 0$ as $r \rightarrow 0.$  For any $g \in \text{VMO},$
$\eta_{g}(r)$ is referred to as the VMO-modulus.  For all conventions regarding VMO we follow \cite{BT}
which in turn follows \cite{MPS}.

\begin{theorem}[Existence of Blowup Limits I] \label{blowup}
Assume $w$ satisfies Equations\refeqn{Always}and\refeqn{Freq}\!\!, and
assume in addition that $a^{ij}$ and $f$ belong to VMO.
Define the usual rescaling
$$w_{\epsilon}(x):=\epsilon^{-2}w(\epsilon x).$$
Then for any sequence $\{\epsilon_m\}\downarrow 0$, there exists a subsequence, a
real number $\mu \in [\lambda^{\ast}, \Lambda^{\ast}],$ and a symmetric matrix $A = (A^{ij})$ with
$$0<\lambda I\leq A\leq \Lambda I$$
such that for all $i,j$ we have
\begin{equation}
\label{eq:vmoconv}
\newavint{B_{\epsilon_m}} a^{ij}(x) dx \rightarrow A^{ij}
\end{equation}
and
\begin{equation}
\label{eq:vmofconv}
\newavint{B_{\epsilon_m}} f(x) dx \rightarrow \mu \;,
\end{equation}
and on any compact set, $w_{\epsilon_m}(x)$ converges strongly in $C^{\alpha}$ and weakly in $W^{1,2}$ to a function $w_{\infty}\in W^{1,2}_{loc}(\R^n)$, which satisfies:
\begin{equation}
\label{eq:blowuplmt}
D_i A^{ij}D_j w_{\infty}= \chisub{ \{ w_{\infty}>0 \} } \mu \ \ \text{on}\ \R^n,
\end{equation}
and has 0 in its free boundary.
\end{theorem}

\begin{pf}  This proof is so similar to the proof of Theorem 6.1 of \cite{BT} that we leave it as an exercise
for the reader.
\end{pf}

\begin{remark}[Nonuniqueness of Blowup Limits]  \label{NBL}
Notice that the theorem does not claim that the blowup limit is unique.
In fact, it is relatively easy to produce
nonuniqueness even in the case with a constant right hand side,
and it was done in \cite{BT} for the nondivergence form case,
but that counter-example can be copied almost exactly
for the divergence form case.  In the case where the coefficients
of $L$ are constant, one can use the counter-example in \cite{B} to
show nonuniqueness of blowup limits when the right hand side is only
assumed to be continuous.
\end{remark}

\begin{theorem}[Caffarelli's Alternative in Measure (Weak Form)]  \label{CAMW}
Assuming again Equations\refeqn{Always}and\refeqn{Freq}\!\!, the limit
\begin{equation}
    \lim_{r \downarrow 0} \frac{ |\Lambda(w) \cap B_r| }{ |B_r| }
\label{eq:densitystatement}
\end{equation}
exists and must be equal to either $0$ or $1/2.$
\end{theorem}

\begin{pf}
Here again our proof is almost identical to the proof of Theorem 6.3 of \cite{BT}, so
we leave it to the reader.
\end{pf}

\begin{definition}[Regular and Singular Free Boundary Points]  \label{RSFBP}
A free boundary point where $\Lambda$ has density equal to $0$ is referred to as \textit{singular}, and
a free boundary point where the density of $\Lambda$ is $1/2$ is referred to as \textit{regular}.
\end{definition}

The theorem above gives us the alternative, but we do not
have any kind of uniformity to our convergence.  Caffarelli
stated his original theorem in a much more quantitative (and
therefore useful) way, and so now we will state and prove a similar
stronger version.  We need the stronger version in order to
show openness and stability under perturbation of the regular
points of the free boundary.

\begin{theorem}[Caffarelli's Alternative in Measure (Strong Form)]  \label{CAMS}
Once again assuming Equations\refeqn{Always}and\refeqn{Freq}\!\!,
for any $\epsilon \in (0, 1/8),$  there exists an $r_0 \in (0,1),$ and a $\tau \in (0,1)$
such that \newline
\indent if there exists a $t \leq r_0$ such that
\begin{equation}
   \frac{|\Lambda(w) \cap B_t|}{|B_t|} \geq \epsilon \;,
\label{eq:bigonce}
\end{equation}
\indent then for all $r \leq \tau t$ we have
\begin{equation}
   \frac{|\Lambda(w) \cap B_r|}{|B_r|} \geq \frac{1}{2} - \epsilon \;,
\label{eq:bigallatime}
\end{equation}
and in particular, $0$ is a regular point according to our definition.
The $r_0$ and the $\tau$ depend on $\epsilon$ and on the $a^{ij},$ but they do \textit{not} depend on the
function $w.$
\end{theorem}

\begin{remark}[Another version]  \label{AnoVer}
The theorem above is equivalent to a version using a modulus of continuity.  In that version there is a universal
modulus of continuity $\sigma$ such that
\begin{equation}
   \frac{|\Lambda(w) \cap B_{\tilde{t}}|}{|B_{\tilde{t}}|} \geq \sigma(\tilde{t})
\label{eq:OtherVersion}
\end{equation}
for any $\tilde{t}$ implies a uniform convergence of the density of $\Lambda(w)$ to $1/2$ once $B_{\tilde{t}}$
is scaled to $B_1.$  (Here we mean uniformly among all appropriate $w$'s.)
\end{remark}

\begin{pf}
Here again we have a proof which is almost identical to the proof of Theorem 6.5 in \cite{BT}.
On the other hand, in an effort to make things more convenient for the reader, since we use this theorem
quite a bit, we will include the proof here.

We start by assuming that we have a $t$ such that Equation\refeqn{bigonce}holds, and by rescaling if necessary,
we can assume that $t = r_0.$  Next, by arguing exactly as in the last theorem, by assuming that $r_0$ is
sufficiently small, and by defining $s_0 := \sqrt{r_0},$ we can assume without loss of generality that
\begin{equation}
\myavinttwo{B_{s_0}} \left| a^{ij}(x) - \delta^{ij} \right| \; dx
\label{eq:atbt}
\end{equation}
is as small as we like.  Now we will follow the argument given for Theorem 4.5 in \cite{B} very closely.

Applying our measure stability theorem on the ball $B_{s_0}$ we have the existence of a function
$u$ which satisfies:
\BVPcbsn{\Delta u = \chisub{\{u > 0\}}\mu}{u \equiv w}{newudef}
and so that
\begin{equation}
|\{\Lambda(u) \Delta \Lambda(w)\} \cap B_{r_0}|
\label{eq:damnsmall}
\end{equation}
is small enough to guarantee that
\begin{equation}
   \frac{|\Lambda(u) \cap B_{r_0}|}{|B_{r_0}|} \geq \frac{\epsilon}{2} \;,
\label{eq:uLbigonce}
\end{equation}
and therefore
\begin{equation}
   m.d.(\Lambda(u) \cap B_{r_0}) \geq C(n) r_0 \epsilon \;.
\label{eq:uLmdbigonce}
\end{equation}
Now if $r_0$ is sufficiently small, then by Caffarelli's $C^{1,\alpha}$ regularity theorem for the obstacle problem
(see \cite{C4} or \cite{C5})
we conclude that $\partial \Lambda(u)$ is $C^{1,\alpha}$ in
an $r_0^2$ neighborhood of the origin. Furthermore, if we rotate coordinates so that
$FB(u) = \{ (x', x_n) \; | \; x_n = g(x') \},$
then we have the following bound (in $B_{r_0^2}$):
  \begin{equation}
  ||g||_{_{C^{1,\alpha}}} \leq \frac{C(n)}{r_0} \;.
  \label{eq:fSmoothBd}
  \end{equation}
On the other hand, because of this bound, there exists a $\gamma < 1$ such
that if $\rho_0 := \gamma r_0 < r_0,$ then
  \begin{equation}
    \frac{ | \Lambda(u) \cap B_{\rho_0} | }{| B_{\rho_0} | }
    \; > \; \frac{1 - \epsilon}{2} \;.
    \label{eq:BrideOfBigLamb}
  \end{equation}
Now by once again requiring $r_0$ to be sufficiently
small, we get
  \begin{equation}
    \frac{ | \Lambda(w) \cap B_{\rho_0} | }{| B_{\rho_0} | }
    \; > \; \frac{1}{2} - \epsilon \;.
    \label{eq:HoundOfBigLamb}
  \end{equation}
(So you may note that here our requirement on the size of $r_0$ will be much smaller than it was before; we need it small
both because of the hypotheses within Caffarelli's regularity theorems and because of the need to shrink the $L^p$ norm of
$|a^{ij} - \delta^{ij}|$ and the $L^1$ norm of $|f - \mu|$ in order to use our measure stability theorem.)

Now since $\frac{1}{2} - \epsilon$ is strictly
greater than $\epsilon,$
we can rescale $B_{\rho_0}$ to a ball with a radius
\textit{close to} $r_0,$ and then repeat.  Since we have
a little margin for error in our rescaling, after
we repeat this process enough times we will have
a small enough radius (which we call $\tau r_0$),
to ensure that for all $r \leq \tau r_0$ we have
$$\frac{ | \Lambda(w) \cap B_r |}{| B_r |}
  \; > \; \frac{1}{2} - \epsilon \;.$$
\end{pf}


\begin{corollary}[The Set of Regular Points Is Open]  \label{TSoRPIO}
Still assuming Equations\refeqn{Always}and\refeqn{Freq}\!\!,
the set of regular points of $FB(w)$ is an open subset of $FB(w).$
\end{corollary}

\noindent
The proof of this corollary is identical to the proof of Corollary 4.8 in \cite{B} except that in place of using Theorem 4.5 of
\cite{B} we use Theorem\refthm{CAMS}from this work.

\begin{theorem}[Existence of Blowup Limits II] \label{blowup2}
We assume Equation\refeqn{Always}\!\!, and we assume $a^{ij}$ and $f$ belong to VMO.
We let
\begin{equation}
    S_r := \{ x \in FB(w) \cap B_r \; : \; x \ \text{is a regular point of} \ FB(w) \; \}
\label{eq:Srdef}
\end{equation}
and we assume $S_{1/2} \ne \phi.$
Let $K \subset \subset S_{1/2},$ let $\{ x_m \} \subset K,$ and let $\epsilon_m \downarrow 0.$

Then there exists a constant $\mu \in [\lambda^{\ast}, \Lambda^{\ast}],$ a constant
symmetric matrix $A = (A^{ij})$ with $0 < \lambda \leq A \leq \Lambda,$
and a strictly increasing sequence of natural numbers $\{ m_j \}$
such that the sequence of functions $\{ w_j \}$ defined by
\begin{equation}
    w_j(x) := \epsilon_{m_j}^{-2}w(x_{m_j} + \epsilon_{m_j}x)
\label{eq:wjdef}
\end{equation}
converges strongly in $C^{\alpha}$ (for some $\alpha > 0$) and
weakly in $W^{1,2}$ on any compact set to a function $w_{\infty}$ 
which satisfies:
\begin{equation}
\label{eq:blowuplmt2}
D_i A^{ij}D_j w_{\infty}= \chisub{ \{ u_{\infty}>0 \} } \mu \ \ \text{on}\ \R^n \;.
\end{equation}
Furthermore $0$ is a regular point of its free boundary.
\end{theorem}

\begin{pf}
The existence of a function $w_{\infty} \geq 0$ satisfying Equation\refeqn{blowuplmt2}and the convergence of
the $w_j$ to $w_{\infty}$ is carried out in exactly the same way as in the proof of Theorem\refthm{blowup}\!\!.
Showing that $0$ is part of the free boundary of $w_{\infty}$ is also proven exactly as in
Theorem\refthm{blowup}\!\!.  It remains to show that $0$
is a regular point of the free boundary.

For the first part, we observe that since each $x_m$ belongs to the regular part of the free
boundary, we know that there exists an $r_m$ such that
\begin{equation}
\label{eq:LamIs3/8}
    \frac{\Lambda(w) \cap B_{r_m}(x_m)}{B_{r_m}} \geq \frac{3}{8} \;.
\end{equation}
There exists a small $\rho > 0$ depending only on the dimension, $n,$ such that if $x \in B_{\rho r_m}(x_m),$ then
\begin{equation}
\label{eq:LamIs1/4}
    \frac{\Lambda(w) \cap B_{r_m}(x)}{B_{r_m}} \geq \frac{1}{4} \;.
\end{equation}

Now the closure of the set $\{ x_m \}$ is compact, and that set is covered by the open balls in the set
$\{ B_{\rho r_m}(x) \}.$ By compactness, the set is still covered by a finite number of these balls, and their
radii have a positive minimum, $\rho_0.$  So, once $\epsilon_{m_j} < \rho_0,$ we know that
\begin{equation}
\label{eq:LamIs1/4Applied}
    \frac{\Lambda(w_j) \cap B_{r}}{B_{r}} \geq \frac{1}{4} \;,
\end{equation}
for all $r$ which are less than $\tau$ times $\rho_0.$ Here $\tau$
is the constant given in the statement of Theorem\refthm{CAMS}\!\!.
From this we can conclude that $0$ must
be a regular point of $FB(w_{\infty}).$
\end{pf}

\begin{remark}[Hausdorff Dimension]   \label{HDim}
Exactly as in \cite{BT}, the arguments above lead to the statement that the free boundary is strongly porous and
therefore has Hausdorff dimension strictly less than $n.$  (See \cite{BT} and see \cite{M} for the definition of
porosity.)
\end{remark}


\newsec{Finer Regularity of the Free Boundary}{SRFB}
In this section we show finer properties of the free boundary at regular points.
Since the counter-examples in \cite{B} and in \cite{BT} are easily extended to the current setting,
we can have regular free boundary points where the blowup limit is not unique.
In spite of this fact, we show that the regular free boundary points enjoy a 
flatness property which is based on Reifenberg flatness.
Reifenberg flatness was introduced by Reifenberg in 
\cite{R}, and is studied in more detail by Toro and 
Kenig in several papers.  (See \cite{KT1} and \cite{KT2} for example.)
For the definitions surrounding Reifenberg vanishing sets we follow
the conventions in section 6 of \cite{B}, but now we must introduce
a notion of sets which are ``relatively Reifenberg flat.''

\begin{definition}[Reifenberg Flatness]  \label{ReifFlat}
Let $S \subset \R^n$ be a locally compact set, and
let $\delta > 0.$  Then $S$ is 
$\mathit{\delta\!-\!\textit{Reifenberg flat}}$ if for each
compact $K \subset \R^n,$ there exists a constant $R_K > 0$
such that for every $x \in K \cap S$ and every $r \in
(0, R_K]$ we have a hyperplane $L(x,r)$ containing $x$
such that
  \begin{equation}
  D_{\mathcal{H}}(L(x,r) \cap B_r(x), \; S \cap B_r(x)) \leq 2r\delta \;.
  \label{DefReif}
  \end{equation}
Here $D_{\mathcal{H}}$ denotes the Hausdorff distance:  If $A, \; B \subset
\R^n,$ then
  \begin{equation}
  D_{\mathcal{H}}(A,B) := \max\{\; \sup_{a \in A} d(a,B) \;, 
  \; \sup_{b \in B} d(b,A) \; \} \;.
  \label{Hdist}
  \end{equation}
We also define
the following quantity, which we call the \textit{modulus of
flatness,} to get a more quantitative and uniform
measure of flatness:
  \begin{equation}
  \theta_K(r) := \sup_{0 < \rho \leq r} \left(
  \sup_{x \in S \cap K} \frac{D_{\mathcal{H}}(L(x,\rho) \cap B_{\rho}(x), \;
  S \cap B_{\rho}(x))}{\rho} \right) \;.
  \label{ThetaFlat}
  \end{equation}
Finally, we will say that $S$ is a \textit{Reifenberg
vanishing} set, if for any compact $K \subset S$
  \begin{equation}
  \lim_{r \rightarrow 0}\theta_K(r) = 0 \;.
  \label{ReifVan}
  \end{equation}
\end{definition}

\begin{definition}[Relatively Reifenberg Flat]  \label{RelReifFlat}
Let $S \subset \R^n$ be a locally compact set, let $K \subset \subset S,$ and
let $\delta > 0.$  Then $K$ is 
\textit{relatively} $\mathit{\delta\!-\!}$\textit{Reifenberg flat with respect to} $\mathit{S}$ if
there exists a constant $R > 0$
such that for every $x \in K$ and every $r \in (0, R]$ we have a hyperplane $L(x,r)$ containing $x$
such that
  \begin{equation}
  D_{\mathcal{H}}(L(x,r) \cap B_r(x), \; S \cap B_r(x)) \leq 2r\delta \;.
  \label{DefRelReif}
  \end{equation}
We also define the \textit{modulus of flatness,} exactly as above,
and then $K$ is \textit{relatively Reifenberg vanishing} if the
modulus of flatness goes to zero as $r$ approaches $0.$
\end{definition}

\begin{remark}  \label{DiffK}
It is worth noting that the compact set $K,$ plays a very different role in the two definitions above.
In the first case, $K$ allows us to look at bounded sets to get uniform bounds on the constant
$R_K$ which bounds the radius, while in the second case, $K$ \textit{is} the set that we want to
show is Reifenberg vanishing, but we are allowing all of $S$ when seeing if we are close to a plane.
As a simple example, a point can never be Reifenberg flat, but viewed as a subset of a plane, it is
relatively $\delta$-Reifenberg flat.
\end{remark}

First we need to show that our measure stability theorem can be used to show uniform closeness of our
solutions to solutions of obstacle problems with constant coefficients and constant right hand side, as long
as we have zoomed in far enough.  In particular, we can say the following:

\begin{theorem}[Uniform Closeness Result]   \label{UCR}
We assume Equation\refeqn{Always}\!\!,
and we let $u \geq 0$ satisfy:
\BVPb{\Delta u = \chisub{\{u > 0\}}\mu}{u \equiv w}{newerudef}
We also assume that there is a fixed constant $\beta,$ and an $\alpha \in (0,1)$
such that $||w||_{C^{\alpha}(\closure{B_1})} \leq \beta.$
For any $\epsilon > 0,$ there exists a $\delta > 0$ such that if
\begin{equation}
   ||a^{ij}(x) - \delta^{ij}||_{L^1(B_1)} < \delta \ \ \ \text{and} \ \ \ ||f(x) - \mu||_{L^1(B_1)} < \delta \;,
\label{eq:propercloseness}
\end{equation}
then
\begin{equation}
   ||w - u||_{L^{\infty}(B_{3/4})} < \epsilon \;.
\label{eq:uniformvictory}
\end{equation}
\end{theorem}

\begin{pf}
Some of the ideas in this proof were inspired by ideas of Li and Vogelius who in turn were following
ideas of Caffarelli.  (See \cite{LV} and \cite{C3}.)
Letting $A(x)$ be the matrix determined by $a^{ij}(x),$ we have in $B_1$
(using ``divergence'' notation):
\begin{alignat*}{1}
  \text{div} &\left[ A(x) \left( \nabla [w(x) - u(x)] \right) \right] \\
         &= f(x) \chisub{ \{ w > 0 \} } - \text{div} \left[ \left( A(x) - I \right) \nabla u(x) \right] - \Delta u \\
         &= f(x) \chisub{ \{ w > 0 \} } - \mu \chisub{ \{ u > 0 \} } -
                             \text{div} \left[ \left( A(x) - I \right) \nabla u(x) \right] \\
         &= f(x) \left( \chisub{ \{ w > 0 \} } - \chisub{ \{ u > 0 \} } \right)
            + \chisub{ \{ u > 0 \} } \left( f(x) - \mu \right)
             + \text{div} \left[ \left( I - A(x) \right) \nabla u(x) \right] \\
         &= I + I\!I + \text{div} \left[ I\!I\!I \right] \;.
\end{alignat*}
After fixing $q \in (n, \infty),$ and by shrinking $\delta$ if necessary, we can use our measure stability theorem
(Theorem\refthm{MeaStab}\!\!) and a simple interpolation, to ensure that the $L^{q/2}$ norm of $I$ on $B_1$
is as small as we like.  Using our assumptions and shrinking $\delta$ if necessary, we can make the $L^{q/2}$
norm of $I\!I$ on $B_1$ as small as we like.  (The fourth line of Equation\refeqn{Always}supplies the
$L^{\infty}$ bound needed for the interpolation.)

To control $I\!I\!I$ we need to shrink the ball slightly.  First we observe that by De Giorgi-Nash-Moser theory
(see Theorem 8.29 of \cite{GT}), there exists an $\alpha^{\prime} \in (0, \alpha)$ such that
\begin{equation}
   ||u||_{C^{\alpha^{\prime}}(\closure{B_1})} \leq C(\beta, \Lambda^{\ast}) \;.
\label{eq:GlobDNM}
\end{equation}
For any fixed $s \in (0,1/16)$ we then have
\begin{equation}
   ||w-u||_{L^{\infty}(\partial B_{1-s})} \leq C(\beta, \Lambda^{\ast})s^{\alpha^{\prime}} \;.
\label{eq:simpleHolder}
\end{equation}
For any $\tilde{q} < \infty,$ we can use Calderon-Zygmund Theory along with the Sobolev Imbedding
to show
\begin{equation}
   ||\nabla u||_{L^{\infty}(B_{1-s})} \leq C(\beta, \Lambda, s) \;.
\label{eq:goodgradbound}
\end{equation}

Considering the boundary value problem that $w - u$ satisfies within $B_{1-s},$ we have the following:
By shrinking $s$ we can make the boundary values as small as we like by Equation\refeqn{simpleHolder}\!\!.
We already have the $L^{q/2}$ norm of $I$ and $I\!I$ as small as we like by making $\delta$ small.
For $I\!I\!I$ we can use Equation\refeqn{goodgradbound}to ensure that
$||\nabla u||_{L^{\infty}(B_{1-s})}$ is under control, and then shrink $\delta$ if necessary to ensure that
$||A - I||_{L^{q}(B_1)}$ is as small as we like.  Applying Theorem 8.16 of \cite{GT} yields the desired
result.
\end{pf}

Now we have a standard corollary for obstacle type problems.

\begin{corollary}[Free Boundaries Are Close]  \label{FBAC}
Assuming Equation\refeqn{Always}again, assuming $u$ is defined as in the previous theorem,
and using $D_{\mathcal{H}}$ as the Hausdorff
distance between sets defined at the beginning of this section, there exists a universal
constant $C$ such that
\begin{equation}
    D_{\mathcal{H}}(FB(w), FB(u)) \leq C\sqrt{\epsilon}
\label{eq:FBsClose}
\end{equation}
where $\epsilon$ is the number given in Equation\refeqn{uniformvictory}\!\!.
\end{corollary}

\begin{pf}
This result is a simple application of the nondegeneracy enjoyed by each function.  Indeed, if there
is a point $x$ where one function is positive and a ball $B_{r}(x)$ where the other function is zero,
then nondegeneracy implies that the max of the first function is $Cr^2$ on $\partial B_{r}(x)$ and this
must be smaller than $\epsilon.$
\end{pf}

Now we prove the main theorem in this paper.

\begin{theorem}[Free Boundary Regularity]  \label{FBR}
Once again we assume Equation\refeqn{Always}and we assume that $a^{ij}$ and $f$ belong to VMO.
As in Equation\refeqn{Srdef}we define $S_r$ to be the set of regular points of the free boundary within
$B_r.$  Let $K \subset \subset S_{1/2}.$  Then
$K$ is relatively Reifenberg vanishing with respect to $S_{1/2}.$
\end{theorem}

\begin{pf}
Fix$\epsilon > 0.$  We will demonstrate that there is a radius $\tilde{r} > 0$ such that for any
$x \in K,$ and any positive $r < \tilde{r}$ there is a hyperplane $H(r,x)$ such that
\begin{equation}
    D_{\mathcal{H}}(FB(w) \cap B_r(x), H(r,x) \cap B_r(x)) \leq r\epsilon \;.
\label{eq:WeWin}
\end{equation}

We start by using the compactness of $K$ in almost the same way as in Theorem\refthm{blowup2}\!\!.
Namely, we know that for every $x \in K$ there exists an $r_x$ such that
\begin{equation}
\label{eq:LamIs3/8rx}
    \frac{\Lambda(w) \cap B_{r_x}(x)}{B_{r_x}} \geq \frac{49}{100} \;.
\end{equation}
Next, there exists a small $\rho > 0$ depending only on the dimension, $n,$ such that if
$y \in B_{\rho r_x}(x) \cap FB(w),$ then
\begin{equation}
\label{eq:LamIs1/4y}
    \frac{\Lambda(w) \cap B_{r_x}(y)}{B_{r_x}} \geq \frac{48}{100} \;.
\end{equation}

Now $K$ is compact, and is therefore covered by the open balls in the set
$\{ B_{\rho r_x}(x) \}.$ By compactness, the set is still covered by a finite number of these balls, and their
radii have a positive minimum, $\rho_0.$  Using Theorem\refthm{CAMS}guarantees that for all $r < \tau \rho_0,$
and for all $x \in K,$ we have
\begin{equation}
\label{eq:wdensg}
    \frac{\Lambda(w) \cap B_{r}(x)}{B_{r}} \geq \frac{48}{100} \;.
\end{equation}
Here $\tau$
is the constant given in the statement of Theorem\refthm{CAMS}\!\!.  Henceforth, the argument becomes
completely independent of whatever point in the free boundary that we wish to consider, so we can fix
$x_0 \in FB(w),$ and show flatness at that point.  Also, given the VMO-modulus $\eta,$ we can be sure
that every quantity that we wish to control below can be shrunk in a uniform and universal way by shrinking
the radius that we are considering.

We consider the situation in $B_{r}(x_0)$ and after a linear
invertible change of coordinates with eigenvalues bounded away from $0$ and $\infty$ in a uniform way
depending only on ellipticity, we can assume that the averages of $a^{ij}$ are $\delta^{ij}$ and the average of
$f$ is $\mu.$  Then we let $u$ solve the boundary value problem:
\BVP{\Delta u = \mu \chisub{ \{ u > 0 \} } }{ u = w }{ourbvp}{ B_{r}(x_0) }
By the $L^1$ closeness of $a^{ij}$ to $\delta^{ij}$ and $f$ to $\mu$ which are controlled by the VMO-modulus
along with our measure stability theorem (Theorem\refthm{MeaStab}\!\!), we can guarantee (by assuming
$r_1$ is sufficiently small) that
\begin{equation}
\label{eq:vdensg}
    \frac{\Lambda(v) \cap B_{r_1}(x_0)}{B_{r_1}} \geq \frac{47}{100} \;.
\end{equation}
Now it follows from Caffarelli's free boundary regularity theorem (see Theorem 7 of \cite{C4} or \cite{C5})
that if $r_2 \leq \tau_2 r_1$ where $\tau_2$ is suitably small, then 
$FB(v) \cap B_{r_2}(x_0)$ is uniformly $C^{1,\alpha}$ in $B_{r_2}(x_0).$  We can also assume
that $FB(v)$ has a free boundary point as close to $x_0$ as we like by using the last corollary (and shrinking $r_1$ again
if needed).  Now zooming in on a uniformly $C^{1,\alpha}$ set will flatten it in a uniform way depending only on how
much one zooms, so after zooming in to $r_3 := \tau_3 r_2,$ where $\tau_3$ will only depend on estimating how
uniformly $C^{1,\alpha}$ functions flatten out as you zoom in, so we can have $FB(v) \cap B_{r_3}(x_0)$ within
$r_3 \cdot \epsilon/2$ of a plane.  Now we invoke Corollary\refthm{FBAC}again to guarantee
that $FB(w)$ is within $r_3 \cdot \epsilon/2$ of $FB(v)$ and we are done.
\end{pf}

\begin{remark}[Choosing $r$]  \label{Chor}
It is worth remarking that the $r_j$ that work for all of the estimates in the last proof must be found \textit{before}
finding the function $u,$ and then in Equation\refeqn{ourbvp}we can use $r = r_3.$
\end{remark}

\begin{remark}[Nondivergence Form Case]   \label{NDFC}
The Theorem above (and the next corollary) can be extended without any difficulty to the nondivergence form
setting.  On the other hand, in the nondivergence form setting, since the functions will have stronger convergence
to their blowup limits, it is very likely that the Weiss-type Monotonicity formula can be used to give an easier proof.
In the divergence form case, the presence of the Dirichlet integral within the Weiss-type monotonicity functional
coupled with the weak convergence in $W^{1,2}$ to the blowup limit makes it difficult to move back and forth from
the original function to its blowup limit.
\end{remark}

\begin{corollary}[Blowup Classification]   \label{BC}
Any blowup found in Theorem\refthm{blowup2}must be homogeneous of degree two, and therefore in the right
coordinate system, it willl be a constant times $(x_n^{+})^2.$
\end{corollary}

\begin{pf}
By Theorem\refthm{FBR}any blowup found in Theorem\refthm{blowup2}will have to be a global solution to the
obstacle problem with a free boundary which is a hyperplane.  Then by applying a combination of the Cauchy-Kowalevski
theorem and Holmgren's uniqueness theorem we conclude (after a possible rotation and change of coordinates)
that the blowup limit is $C(x_n^{+})^2.$
\end{pf}


\bibliographystyle{abbrv}
\bibliography{simple}

\end{document}